\documentclass[12pt,a4paper]{article}
\usepackage{graphicx}
\usepackage[english]{babel}
\usepackage{amsmath}
\usepackage[a4paper, total={6in, 9in}]{geometry}
\usepackage{amsfonts}
\usepackage{amssymb}
\usepackage{amsthm}
\usepackage{amstext} 
\usepackage{array}
\newtheorem{lemma}{Lemma}
\usepackage{pdflscape}
\usepackage{rotating}
\usepackage{float}
\usepackage{hyperref}

\newtheorem*{theorem*}{Theorem}

\newtheorem{defi}{Definition}
\newtheorem{prop}{Proposition}
\newtheorem{coro}{Corollary}
\newtheorem{statement}{Statement}
\newtheorem{notat}{Notation}
\newtheorem{conj}{Conjecture}
\newcolumntype{C}{>{$}c<{$}}
\title{On axial algebras with $3$ eigenvalues}
\author{Vsevolod A. Afanasev}
\date{}

\begin{document}
\maketitle

\begin{abstract}
We study axial algebras, that is, commutative non-associative algebras generated by idempotents whose adjoint actions are semisimple and obey a fusion law. Considering the case, when said adjoint actions having $3$ eigenvalues and the fusion law is the least restrictive, we describe $2$-generated and $3$-generated algebras and prove some of their general properties.
\end{abstract}

\let\thefootnote\relax\footnotetext{The work was supported by the Theoretical Physics and Mathematics Advancement Foundation ”BASIS”.}
The theory of axial algebras has been introduced in 2015 in the work of Hall, Rehren and Shpectorov \cite{Axial}. These algebras are commutative, non-associative and are generated by \textit{axes} -- idempotents, such that for any axis $a$ its multiplication operator $ad_a$ splits the algebra into the direct sum of eigenspaces that multiply according to the certain \textit{fusion law} (see precise definitions below). The algebras were introduced as objects, axiomatising key properties of vertex operator algebras, but have since been investigated from multiple points of view. As of now, most of the work done in the area was motivated by two major examples: the class of algebras of Jordan type, which includes idempotent-generated Jordan algebras and the algebras of Monster type, among them the famous Griess algebra. Many insightful results were obtained for both cases, which include general properties as well as classifications for specific families of algebras  (for more details see survey \cite{ShpecMcSurv}).

One of the more significant differences between the two examples is the fact that $ad_a$ has $3$ eigenvalues for the algebras of Jordan type and $4$ eigenvalues for the algebras of Monster type. In fact, algebras of Jordan type can be viewed as algebras of Monster type with one of the eigenvalues multiplying trivially in the fusion law. However greater generality is compensated by greater difficulty of studying this class. Recently, an application of axial framework to the idempotent-generated pseudo-composition algebras and the train algebras of rank 3 allowed to unite these families as algebras of pseudo-composition type \cite{GorMamStar}. This class also has $3$ eigenvalues for $ad_a$ but differs from the class of algebras of Jordan type. Nonetheless, some techniques applied to the study of the Jordan type worked for the pseudo-composition type as well. Therefore it is natural to consider the $3$-eigenvalue case in the greatest possible generality in order to further study the properties discovered in the examples and start a classification of finitely generated algebras, hoping to establish techniques that may aid in future research in the field of axial algebras and yield more interesting examples.

In this paper we study the fusion law $\mathcal{J}(\alpha,\beta)$ (see Table \ref{Jab}), which is the most general fusion law for the $3$ eigenvalue case. It is «the most general» in the sense that we employ only the most natural assumptions, namely that the fusion law is $\mathbb{Z}/2\mathbb{Z}$-graded and that $1*\lambda=\lambda$ for any $\lambda$ from the fusion law (see the precise definitions below). We will answer some natural questions about the algebras, satisfying this fusion law and provide classification results for $2$-generated and $3$-generated such algebras. The fusion law we are considering (and the results we obtain for it) generalises the fusion laws for the algebras of Jordan and pseudo-composition types in particular.
Similar work for the case of $2$-generated algebras has been done earlier and is a part of Whybrow's unpublished work \cite{Whybrow}, which, as we will show later, contains some errors.

The structure of the paper goes as follows: we list most of our definitions and necessary notation in Section 1. Then in Section 2 we classify the $2$-generated $\mathcal{J}(\alpha, \beta)$-algebras and state and prove some general results about certain classes of such algebras. Finally, Section $3$ contains the study and classification (including the construction of a spanning set and the algorithm to construct a multiplication table) of $3$-generated $\mathcal{J}(\alpha,\beta)$-algebras, that satisfy the $(\star)$ condition (stated in Section 2). 

\paragraph*{Acknowledgements.}
The author would like to thank Alexey Staroletov for reading the initial draft of the paper and pointing out multiple inaccuracies and for all of his incredibly useful remarks and suggestions.

\section{Notation and definitions}

Here we will provide some definitions that are standard to the theory of axial algebras.

\begin{defi}
Let $\mathbb{F}$ be a field. Then the fusion law $\mathcal{F}$ is a subset of $\mathbb{F}$ with a map $*: \mathcal{F} \times \mathcal{F} \to 2^\mathcal{F}$. 
\end{defi}

Before we introduce the definition of an axis we need to fix some notation: first for an element $a$ of an $\mathbb{F}$-algebra $A$ by $ad_a$ we denote the left multiplication operator i.e. $ad_a(x)=ax$. By $A_\lambda(a)$ we denote the (possibly trivial) $ad_a$-eigenspace, corresponding to $\lambda$. And finally by $A_S(a)$, where $S \subseteq \mathbb{F}$ we denote a direct sum $\bigoplus_{\lambda \in S} A_\lambda(a)$.

\begin{defi}
    Let $A$ be a commutative non-associative $\mathbb{F}$-algebra and let $\mathcal{F}$ be a fusion law. Then an idempotent $a \in A$ is called an $\mathcal{F}$-axis if
    \begin{enumerate}
        \item With respect to the map $ad_a$, the algebra can be split into a direct sum of eigenspaces as $A=\bigoplus\limits_{\lambda \in \mathcal{F}}A_\lambda(a).$
        \item For any $\lambda, \mu  \in \mathcal{F}$, we have that $A_\lambda(a) A_\mu(a) \subseteq A_{\lambda * \mu}(a)$.
    \end{enumerate}
\end{defi}

By definition, $1$ is always an $ad_a$-eigenvalue for an axis $a$. We call an axis \textit{primitive} if $A_1(a)=\langle a \rangle $. 

Here and throughout the rest of this paper we use $\langle a_1,\ldots, a_n \rangle$ to denote a vector space, generated by $a_1,\ldots, a_n$. 

\begin{defi}
    We call an algebra $A$ over field $\mathbb{F}$ a (primitive) axial algebra for the fusion law $\mathcal{F}$ if it is generated by a set of (primitive) $\mathcal{F}$-axes.
\end{defi}

This definition depends on the generating set. We will provide such sets explicitly where necessary.  

Most of the time it is convenient to present the fusion law as a table. In particular, the fusion law that we consider in this paper can be written as follows:

\[
\label{Jab}
\mathcal{J}(\alpha,\beta) = \begin{tabular}{|C|C|C|C|}\hline
   * & 1 & \alpha & \beta \\
   \hline \hline
    1 & 1 & \alpha& \beta \\
  \hline
  \alpha & \alpha & 1,\alpha & \beta \\
  \hline
  \beta & \beta &  \beta & 1,\alpha \\
  \hline
  \end{tabular}\ \ \ \ \ \ \ \ 
\]
\begin{center}
    \text{Table 1: The fusion law } $\mathcal{J}(\alpha,\beta)$
\end{center}

For convenience, we omit the set brackets in each cell of this table.

We assume that all eigenvalues $1,\alpha,\beta$ are distinct.
We call the axial algebra, satisfying this fusion law a $\mathcal{J}(\alpha,\beta)$-algebra.

As was stated earlier, this fusion law is a generalisation of both the $\mathcal{J}(\eta)$ fusion law and $\mathcal{PC}(\eta)$ fusion law, which are the fusion laws for algebras of Jordan type and of pseudo-composition type respectively.

One can note that $\mathcal{J}(\alpha,\beta)$ (and as a result, any axial algebra with that fusion law) can be $\mathbb{Z}/2\mathbb{Z}$-graded with $A_1(a) \oplus A_\alpha(a)$ and $A_\beta(a)$ forming two parts of a grading for any axis $a$. In other words, any element $x$ of $A$ can be written as $x_+ + x_-$ with respect to $a$. We note that all such fusion laws with at most $3$ eigenvalues were classified in \cite{Whybrow}. This grading allows for a construction of an automorphism $\tau_a$ such that
$$\tau_a(x)=x_+- x_-.$$
\begin{defi}
    Let $A$ be an $\mathcal{J}(\alpha,\beta)$-axial algebra. Write $x \in A$ as $x_+ + x_-$, where $x_+ \in A_1(a) \oplus A_\alpha(a)$ and $x_- \in A_\beta (a)$ with respect to axis $a$. Then the map $\tau_a: x_+ + x_- \to x_+ - x_-$ is called a Miyamoto involution, associated to $a$.
\end{defi}
It is clear from the definition that $\tau_a^2=id$ for any axis $a$. Moreover, an important property of Miyamoto involutions is that they take axes to axes with $\tau_a$ sending $A_q(b)$ to $A_q(\tau_a(b))$. 
\section{Algebras generated by 2 axes}

We will first consider the case of a $\mathcal{J}(\alpha,\beta)$-algebra $A$ being generated by primitive axes $a,b$ and start by providing a classification of such algebras. More specifically:

\begin{prop}\label{2gend}
    Assume $A$ is a $\mathcal{J}(\alpha,\beta)$-axial algebra generated by primitive axes $a,b$. Then $a$ is spanned by $a,\ b,\ \tau_a(b)$ or equivalently by $a,\ b,\ ab$. Moreover, if $b=\delta a+b_\alpha+b_\beta$ where $\delta \in \mathbb{F}$, $b_\alpha \in A_\alpha(a)$, $b_\beta \in A_\beta(a)$, then 
    $A_\alpha(a)=\langle (\beta-1) \delta a-\beta b+ab \rangle,\ A_\beta(a)=\langle (\alpha-1) \delta a-\alpha b + ab \rangle$.
\end{prop}

\textit{Proof}. First we rewrite one of the axes with respect to another.
$$b=\delta a+b_\alpha+b_\beta,$$
where $\delta  \in F$. 
Then
$$ab= \delta a+\alpha b_\alpha+ \beta b_\beta.$$

We use these to derive
$$b_\alpha=\frac{ \delta (\beta-1)a-\beta b+ ab}{(\alpha-\beta)},\ b_\beta=\frac{-\delta (\alpha-1)a+\alpha b-ab}{(\alpha-\beta)}.$$

From these we compute 
$$a(ab)=a(\delta a+\alpha b_\alpha + \beta b_\beta)= \delta a+\alpha^2 b_\alpha + \beta^2 b_\beta = (\alpha-1)(\beta-1) \delta a-\alpha \beta b+(\alpha+\beta)ab.$$
Symmetrically, we can rewrite $a$ with respect to $b$ as
$$a=\theta b+a_\alpha+a_\beta.$$

Repeating the above operations we can derive that
$$b(ab)=(\alpha-1)(\beta-1)\theta b-\alpha\beta a+(\alpha+\beta)ab.$$

Now consider $$\tau_a(b)=\delta a+b_\alpha-b_\beta=b-2b_\beta=\frac{2\delta (\alpha-1)a-(\alpha+\beta)b+2ab}{(\alpha-\beta)}.$$ 

Since $\tau_a(b)$ is also an axis, we can utilise the identity $\tau_a(b)^2=\tau_a(b)$, 
which can be rewritten as follows:
$$\frac{(2\delta (\alpha-1)a-(\alpha+\beta)b+2ab)^2}{(\alpha-\beta)^2}-\frac{2\delta (\alpha-1)a-(\alpha+\beta)b+2ab}{(\alpha-\beta)}=0.$$

We now aim to present $4(ab)^2$ via $a,b$ and $ab$. 
For this we need to use the expressions of $a(ab)$ and $b(ab)$,
which include $\delta$ and $\theta$ respectively. 
More specifically
$$4(ab)^2=\left(2\delta (\alpha-\beta)(\alpha-1)-4(\alpha - 1)^2(2\beta - 1)\delta ^2 - 4\alpha \beta(\alpha + \beta)\right)a+$$
$$+\left(-(\alpha + \beta)(\alpha-\beta)-\Bigl((\alpha+\beta)^2-8\delta (\alpha-1)(\alpha\beta)-4(\alpha+\beta)(\alpha-1)(\beta-1)\theta \Bigr)\right)b+$$
$$+\left(2(\alpha-\beta)-4 \delta(\alpha - 1) (\alpha + \beta) + 4 (\alpha + \beta)^2\right)ab.$$

After this we can compute the spectrum and eigenvectors of $ad_a$ to get that
$$v_\alpha = (\beta-1)\delta a-\beta b+ab,\ v_\beta=(\alpha-1)\delta a-\alpha b + ab.$$
The eigenvectors for $ad_b$ can be written symmetrically as $v^\prime_\alpha=(\beta-1)\theta b-\beta a+ab,\ v^\prime_\beta=(\alpha-1)\theta b-\alpha a+ ab$. $\square$

The next step in studying a class of axial algebras is traditionally an attempt to construct a Frobenius form.

\begin{defi}
    Suppose $A$ is an axial algebra for the fusion law $\mathcal{F}$. Then a non-zero symmetric bilinear form $(\cdot ,\cdot)$ on $A$ is called a Frobenius form if for any $u,v,w \in A$ we have that
    $(u,vw)=(uv,w)$.
\end{defi}

In the case of algebras of both Jordan and pseudo-composition types, the form was constructed due to the fact that the projections of axes onto each other's $1$-eigenspaces were equal. Here we show that this is not always the case for the $\mathcal{J}(\alpha,\beta)$ fusion law.

\begin{prop}\label{xy}
Let $A$ be a primitive $\mathcal{J}(\alpha, \beta)$-axial algebra and let $a,b \in A$ be the linearly independent primitive axes. Write $b=\delta a+ b_\alpha+ b_\beta$, where $b_\alpha \in A_\alpha(a),\ b_\beta \in A_\beta(a)$. 
Then one of the following statements holds:
\begin{enumerate}
    \item[(i)] The projection of $a$ on $A_1(b)$ equals to $\delta b$.
    \item[(ii)] $A$ is a 2-dimensional algebra, such that $A_\alpha(a)=A_\beta(b)=0$, $ab=\alpha a + \beta b$ and the condition $2\alpha \beta=-(\beta-1)$ holds for its fusion law eigenvalues.
\end{enumerate}
These options are mutually exclusive.
\end{prop}
\textit{Proof.}
Assume first that $a,\ b,\ ab$ are linearly independent. Then consider the expression for $4(ab)^2$ above. Note that it was written out via the identity for $\tau_a(b)$. 
If we consider $\tau_b(a)$ instead we will get a symmetric result with $\delta$ and $\theta$ being replaced. 
Equating the coefficients of $ab$ for these presentations gives us that
$$4(\alpha-1)(\delta-\theta)(\alpha+\beta)=0,$$
meaning that either $\delta=\theta$ or $\alpha=-\beta$.
    
Assume that $\alpha=-\beta$. We have that, in particular, $a(ab)$ and $b(ab)$ lie in $\langle a,b \rangle$, 
more specifically
$$a(ab)=\delta a-\alpha^2(\delta a-b),\ b(ab)=\theta b-\alpha^2(\theta b-a).$$

We express $b(ab)$ in two ways, namely 
$$b(ab)=(\delta a+b_\alpha+b_{-\alpha})(\delta a+\alpha b_\alpha -\alpha b_{-\alpha})=\delta ^2a+\delta \alpha b_\alpha +\alpha^2 \delta  b_\alpha +\alpha b_\alpha^2-\alpha b^2_{-\alpha}+$$
$$\delta \alpha^2 b_{-\alpha}-x\alpha b_{-\alpha}.$$

The "even"\ and "odd"\  parts of the equation grading-wise are displayed on the first and the second line respectively.

Now we transform the previously displayed equation:
$$b(ab)=(\theta-\alpha^2 \theta) (\delta a+b_\alpha +b_{-\alpha})-\alpha^2a.$$

Equating the odd parts gives us
$$(\delta \alpha(\alpha-1)+\theta (\alpha^2-1))b_{-\alpha}=0.$$
Performing a similar argument for $a(ab)$ we obtain
$$(\theta \alpha(\alpha-1)+\delta (\alpha^2-1))a_{-\alpha}=0.$$
Since $a,\ b,\ ab$ are linearly independent, $b_{-\alpha} \neq 0$ and $a_{-\alpha} \neq 0$, therefore we get that $\delta=\theta$.

If we assume that $a,b$ and $ab$ are linearly dependent, then $A_t(a)=0$, where $t \in \{\alpha,\beta\}$, implying that either $v_\alpha$ or $v_\beta$ is equal to $0$. Symmetrically, $A_q(b)=0$ for some $q \in \{\alpha,\beta\}$, meaning that one of the corresponding eigenvectors $v^\prime_\alpha$ or $v^\prime_\beta$ equals to $0$. If $t=q$ then by symmetry of eigenvectors we get that $\delta=\theta$. If $t \neq q$, assume, without loss of generality that $t=\alpha$. Then note that $ab=-(\beta-1)\delta a+\beta b$. Since $A_\alpha(a)=0$ we have that $\beta * \beta = 1$ and so $v_\beta^2=(\alpha-\beta)^2(\delta a-b)^2$ must lie in $\langle a \rangle$, which can happen only if $-2\beta \delta +1=0$. On the other hand, since $A_\beta(b)=0$, $ab=\alpha a-(\alpha-1)\theta b$, therefore $(-(\beta-1)\delta -\alpha)a+(\beta+(\alpha-1)\theta)b=0$, implying in turn that 
$$\delta=\frac{-\alpha}{\beta-1}, \theta=\frac{-\beta}{\alpha-1}.$$
Which, combined with the previous identity gives us that $2\alpha \beta =-(\beta-1)$.
This sets $\delta=\frac{1}{2}(2\alpha+1)$, while $\theta=-\frac{1}{(\alpha-1)(2\alpha+1)}$, leading to $ab=\alpha a+\beta b$. The fusion laws for $b$ are then automatically satisfied, while $v_\beta^2 \in \langle a \rangle$ as shown before, so the fusion law is satisfied for $a$ as well.

If we instead assume that $A_\beta(a)=A_\alpha(b)=0$ we get that $ab=\alpha b+ \beta a$.

If we were to assume that $\delta=\theta$ in this algebra, then, equating the values of $\delta$ and $\theta$ as shown above we get a system of equations
$$\begin{cases}
 \alpha(\alpha-1)=\beta(\beta-1)\\
2\alpha \beta=-(\beta-1)
    \end{cases},
$$
which has solutions $(-1,-1),\ (0,1),\ (\frac{1}{2},\frac{1}{2})$, none of which fit our requirements, thus proving the mutual exclusiveness of the statements in the proposition. $\square$

This shows us that there exist $2$-dimensional, two-generated algebras, on which there is no Frobenius form, whose radical does not contain axes. Due to this we will exclude this family from further consideration in this paper by introducing the following notation.

\begin{notat}
    We say that a primitive $\mathcal{J}(\alpha,\beta)$-axial algebra $A$ satisfies the $(\star)$ condition if for any two primitive axes $a,\  b \in A$, the projection of $a$ onto $A_1(b)$ is equal to $\delta b$ if $b=\delta a+b_\alpha+b_\beta$, $b_\alpha \in A_\alpha(a),\ b_\beta \in A_\beta(a)$.
\end{notat}

Per Proposition \ref{xy}, this condition does not hold only for certain representatives of a family of algebras with a specific condition on the eigenvalues.
The result of Proposition \ref{xy} also contradicts the results of \cite{Whybrow}, where it is claimed that $\delta=\theta$ always.

\begin{coro}\label{xvals}
Assume that $a,\ b$ are primitive $\mathcal{J}(\alpha,\beta)$ axes with the set $a,\ b,\ ab$ being linearly independent. Then we have that at least one of the following is true:
\begin{enumerate}
    \item $\beta=\frac{1}{2},$
    \item $\delta=-\frac{\alpha+\beta}{2(\alpha-1)}.$
\end{enumerate}
\end{coro}

\textit{Proof.}
Computing $(ab)^2$ via $\tau_a(b)$ and $\tau_b(a)$ and subtracting the corresponding coefficients for $a$ yields the identity 
$$(2\beta-1)(x(\alpha-1)+\alpha)(\delta(\alpha-1)+\frac{\alpha+\beta}{2})=0.$$

On the other hand, the determinant of the matrix $\left[a,v_\alpha,\ v^2_\alpha\right]$, with the last vector equal to
$$v_\alpha^2=(\beta-1)^2\delta^2a+\beta^2b+(ab)^2+2\bigl(-\beta(\beta-1)\delta ab+ (\beta-1)\delta a(ab)-\beta b(ab)\bigr)$$

has the determinant equal to $-\frac{1}{2}(2\beta-1)(\alpha-\beta)(2(\alpha-1)\delta+(\alpha+\beta)),$
which is zero iff $\beta=\frac{1}{2}$ or $\delta=\frac{\alpha+\beta}{2(\alpha-1)}$.  

Computing the determinant of the matrix $\left[a,v_\alpha,v^2_\beta \right]$ yields the same restrictions.
Thus the fusion law is satisfied when either of these equalities takes place, proving the result. $\square$

The value $\delta=-\frac{\alpha}{\alpha-1}$ will be making a comeback in the 2-dimensional cases. However, we need to establish the existence of the Frobenius form, where possible, beforehand.
\begin{prop}\label{frobform}
Let $A$ be a primitive $\mathcal{J}(\alpha,\beta)$ axial algebra. Then the following is true.
\begin{enumerate}
    \item $A$ is spanned (as a vector space) by its primitive axes;
    \item If $A$ satisfies the $(\star)$ condition, then $A$ admits a unique Frobenius form such that $(a,a)=1$ for any primitive axis $a$. More specifically, if $b=\delta a+ b_\alpha + b_\beta$, then $(a,b)=\delta$.
\end{enumerate}
\end{prop}

\textit{Proof.}
Given Proposition \ref{xy}, these results can be proved by using the classical arguments which can be found in for instance \cite[Proposition 3,4]{GorMamStar} $\square$.

Proposition \ref{frobform}, together with Corollary \ref{xvals}, show us that outside of $\beta=\frac{1}{2}$ the values of the Frobenius form between axes are heavily restricted. This generalises a well-known result on Frobenius form values for Matsuo algebras. Said result together with the fact that some algebras of Jordan type half are  Matsuo algebras shows that the statements in Corollary \ref{xvals} are not mutually exclusive.

\subsection*{Two-dimensional algebras.}

Here we study the conditions, under which the $2$-generated $\mathcal{J}(\alpha,\beta)$-algebras that satisfy the $(\star)$ condition are $2$-dimensional. That is, we shall assume that $a$ and $b$ are distinct linearly independent axes, while the set $\{a,\ b,\  ab\}$ is linearly dependent. Because of this we have that one of either $b_\alpha, b_\beta$ and one of $a_\alpha=0,\ a_\beta$ equals to $0$.

Under these restrictions the case $b_\alpha=a_\beta=0$ has already been considered in the proof of Proposition 2, yielding no algebras satisfying the $(\star)$ condition.

Therefore assume $b_\alpha=a_\alpha=0$. We then have that $\delta(\beta-1)a-\beta b+ ab=0;\; \delta(\beta-1)b-\beta a+ ab=0$.

We then have that $\delta(\beta-1)a-\beta b-\delta(\beta-1)+\beta a=0$, implying that $\delta(\beta-1)+\beta=0$ (the coefficients for $a$ and $b$ only differ by a sign). In other words $\delta=\frac{-\beta}{\beta-1}$.

Substituting this into $b_\alpha$ we get that $ab=\beta(a+b)$.

Finally, for $b_\beta=a_\beta=0$ we have, through similar reasoning that $ab=\alpha(a+b)$.

Conversely, suppose that $A$ is a $2$-dimensional commutative algebra generated by 2 idempotents $a$ and $b$.

We first assume that $ab=\alpha(a+b)$. Then the $\alpha a + (\alpha-1) b$ is an $\alpha$-eigenvector of adjoint of $a$. Similarly $\alpha b + (\alpha-1)a$ is an $\alpha$-eigenvector for the adjoint of $b$. Then as we have $A=A_1(a)\oplus A_\alpha(a)=A_1(b)\oplus A_\alpha(b)$ and $A_\beta(a)=0= A_\beta(b)$, the fusion rules are trivially satisfied for both $a$ and $b$.

Finally, assume $ab=\beta(a+b)$. Then $\beta a+ (\beta-1) b$ is a $\beta$-eigenvector of adjoint of $a$. Similarly, $\beta b + (\beta-1)a$ is a $\beta$ eigenvector for $b$.
In other words $A=A_1(a)\oplus A_\beta(a)=A_1(b)\oplus A_\beta(b)$. As $A_\alpha(a)=0$ we have that $A_\beta(a)A_\beta(a) \subseteq A_1(a)$ or, in other words $(\beta a+ (\beta-1)b)^2 \in \langle a \rangle$. We have that $(\beta a+ (\beta-1)b)^2=\beta^2 a+ 2\beta^2(\beta-1)(a+b)+ (\beta-1)^2 b$. Since $a$ and $b$ are linearly independent, we have that $2\beta^2(\beta-1)+(\beta-1)^2=0$, or, if we divide by $\beta-1$, which is not $0$, we have that $2 \beta^2=-(\beta-1)$, which implies that $\beta=-1$ or $\beta=\frac{1}{2}$. Due to symmetry, same relations can be obtained from the decomposition of $A$ with respect to the adjoint operator of $b$.

We have arrived at the following

\begin{prop}\label{2dim}
    
Assume that $A$ is a $2$-dimensional axial $\mathcal{J}(\alpha, \beta)$-algebra satisfying the $(\star)$ condition, that is generated by primitive axes $a$ and $b$. Then one of the following is true:
\begin{itemize}
    \item $ab=\alpha(a+b)$. Here $\delta=\frac{-\alpha}{\alpha-1}$,
    \item $ab=\beta(a+b)$ and $\beta=-1, \frac{1}{2}$. Here $\delta=\frac{-\beta}{\beta-1}$.
\end{itemize}
\end{prop}

This shows that for any $\alpha$ and $\beta$ there is at least one family of $2$-dimensional axial algebras for a given fusion law.

Note that this result partially differs from the result obtained for the fusion law $\mathcal{PC}(\eta)$ \cite[Proposition 2]{GorMamStar}, which is stricter than the one we are employing and therefore forbids the existence of some algebras.

\begin{coro}
Assume that $A$ is a $\mathcal{J}(\alpha,\beta)$-axial algebra satisfying the $(\star)$ condition and for any primitive axis $a \in A$  $A_\alpha(a)^2 \subseteq A_1(a)$. Then the first family of algebras in the above Proposition exists only for $\alpha=-1, \frac{1}{2}$.    
\end{coro}

\textit{Proof.}
If we assume $A_\alpha^2 \subseteq A_1(a)$ then, since the $\alpha$-eigenvectors are obtained from $\beta$-eigenvectors by replacing $\beta$ with $\alpha$, an argument, similar to the one used in the previous proof can be utilised to show that $\alpha^2=-\frac{\alpha-1}{2}$. $\square$

\subsection*{Baric algebras}

Another reasonable inquiry is to check, when the $\mathcal{J}(\alpha,\beta)$ algebras can be baric.

Recall that the axial algebra $A$ is called \textit{baric} if there exists a homomorphism $w: A \to \mathbb{F}$, such that $w(a)=1$ for any primitive axis $a \in A$.

In the case of $\mathcal{J}(\alpha,\beta)$-algebras the following statement holds:

\begin{prop}
    Let $A$ be the primitive $\mathcal{J}(\alpha,\beta)$-axial algebra satisfying the $(\star)$ condition. Then $A$ is baric iff $(a,b)=1$ for any primitive axes $a,b$.
\end{prop} 

\textit{Proof.} 
First assume that $A$ is baric with $w:A \to \mathbb{F}$ being the $\mathbb{F}$-algebra homomorphism from the definition. Then note that the form $(a,b)=w(a)w(b)$ is the Frobenius form, which is unique by Proposition \ref{frobform}. In particular $(a,b)=w(a)w(b)=1$.
Now assuming that $(a,b)=1$ for any primitive axis, we can again utilise Proposition \ref{frobform} to check that for any primitive axis $a$ $w_a(x)=(a,x)$ is an $\mathbb{F}$-algebra homomorphism, that is $w_a(x)=w_b(x)$ for any primitive axis $b$ and $(x,y)=w_a(xy)$. Indeed, assume that $x_1,\ldots x_n$ is a basis of $A$ consisting of primitive axes. Then its clear that $w_a(x)=(a,x)=(b,x)=w_b(x)$, since $x$ can be presented as a linear combination of primitive axes $x_1,\ldots x_n$ and by our assumption $(a,x_i)=(b,x_i)=1$ for all $i=1\ldots n$. To prove the second claim we expand 
$$w_a(xy)=(a,xy)=(a,(\delta_1x_1+\ldots\delta_n x_n)(\theta_1x_1+\ldots \theta_n x_n)).$$ 
This is clearly a sum of elements of form $\delta_i \theta_j (a,x_ix_j),\ 1\leq i,j \leq n$. Noting that by previous statement $(a,x_ix_j)=(x_i,x_ix_j)$ we get that this product is equal to 
$$\sum_{1\leq i,j \leq n}\delta_i\theta_j,$$ and this can be easily confirmed to be equal to $(\delta_1x_1+\ldots\delta_n x_n,\theta_1x_1+\ldots \theta_n x_n)$ as all involved axes are primitive.  $\square$


This shows us that in order to find the baric algebras it is possible to consider the subalgebras $\langle \langle a, b \rangle \rangle$, since if $A$ is a baric algebra, then any such subalgebra should be baric as well.

Initially we will compute $(\alpha-\beta)^2 b_\alpha^2$. For convenience we will denote $(a,b)$ by $\delta$. 
After necessary computations we get
$$(\alpha-\beta)^2 b_\alpha^2 = (\alpha-\beta)((\delta^2(-2\alpha \beta+\alpha+\beta)-\alpha\beta+\frac{1}{2}\delta(\alpha -1))a+(\delta(3\alpha \beta-\alpha-\beta+1)-\frac{1}{2}\alpha -\beta)b-$$
$$-(\alpha \delta-\alpha -\beta +\delta-\frac{1}{2})ab).$$

By the fusion rules, this product lies in $A_1(a)\oplus A_\alpha(a)$. Since the first term of the sum obviously lies in $A_1(a)$ we need to unpack the remaining two terms. In particular, the odd parts of these elements sum to $0$.

More specifically
$$(\delta(3\alpha \beta-\alpha-\beta+1)-\frac{1}{2}\alpha -\beta-\beta(\alpha \delta-\alpha -\beta +\delta-\frac{1}{2}))b_\beta=0.$$

This can be transformed into 
$$\frac{1}{2}(\beta-\frac{1}{2})(\alpha(1+2\delta)+\beta-2\delta)b_\beta=0.$$

First we consider the case $b_\beta=0,$ having $ab=-\delta(\alpha-1)a+\alpha b$. This would imply that our algebra is $2$-dimensional and has $a_\beta=0$ (per Proposition \ref{2dim}). This implies that $\delta=\frac{\alpha}{\alpha-1}$ As we want the algebra to be baric, we can then put $\delta=1$, yielding $\alpha=\frac{1}{2}$.

In the other case, since we want our algebras to be baric, we may again assume that $\delta=1$, obtaining the following

\begin{prop}
    
If $A$ is a baric $\mathcal{J}(\alpha,\beta)$ algebra satisfying the $(\star)$ condition, then one of the following is true:
\begin{itemize}
    \item $\alpha=\frac{1}{2},$
    \item $\beta=\frac{1}{2},$
    \item $\beta=2-3\alpha.$
\end{itemize}
\end{prop}

A similar argument provides conditions for existence of the so-called flat algebras, i.e. the algebras such that $(a,b)=0$ for any $2$ primitive axes. These appear as certain exceptions in the works concerning the so-called solid subalgebras in algebras of Jordan type.

\begin{coro}
    If $A$ is a flat $\mathcal{J}(\alpha,\beta)$ algebra, satisfying the $(\star)$ condition then one of the following is true:
    $\beta=\frac{1}{2},\ \alpha=0$ or $\alpha=-\beta$.
\end{coro}

As was mentioned before, the paramount examples of $J(\alpha,\beta)$-algebras are the algebras of Jordan type and of the pseudo-composition type. These cases are different for a variety of reasons and we will show that this difference prevents the techniques used for their classification from being applied to the general case.
We will start with a technical lemma, as we need to calculate the values of the Frobenius form between projections of arbitrary elements.

Let $d=\gamma a+d_\alpha+d_\beta$, $f=\varepsilon a+ f_\alpha+f_\beta$.

\begin{lemma}
For the $\mathcal{J}(\alpha, \beta)$-algebra $A$, that satisfies the $(\star)$ condition, let $a$ be a primitive axis and $d,f \in A$. Then we have that
\begin{itemize}
\item $(d_\alpha , f_\alpha)=\frac{\gamma \varepsilon (\beta-1) - \beta(d,f)+(a,df)}{\alpha - \beta},$ 
\item $(d_\beta, f_\beta)=-\frac{\gamma \varepsilon(\alpha-1)-\alpha(d,f)+(a,df)}{\alpha-\beta}.$
\end{itemize}
\end{lemma}

\textit{Proof.} As before, we can rewrite to show that 

$$d_\alpha = \frac{\gamma (\beta-1)a-\beta d+ad}{(\alpha-\beta)},\ f_\alpha = \frac{\varepsilon (\beta-1)a-\beta f+af}{(\alpha-\beta)}.$$

In other words $$(\alpha-\beta)^2(d_\alpha, f_\alpha)=(\gamma (\beta-1)a-\beta d+ad,\varepsilon (\beta-1)a-\beta f+af).$$
After gathering the terms we get
$$\gamma \varepsilon (\beta-1)(\beta-1-\beta+1+1-\beta)+\beta^2(d,f)-2\beta(a,df)+(ad,af)=$$
$$=-\gamma\varepsilon(\beta-1)^2+\beta^2(d,f)-2\beta(a,df)+(ad,af).$$

We can compute $(ad,af)=(a(ad),f)$. This is equal to 
$$((\alpha-1)(\beta-1)\gamma a-\alpha \beta d + (\alpha+\beta) ad,f)=(\alpha-1)(\beta-1)\gamma \varepsilon - \alpha \beta (d,f) + (\alpha+\beta) (a,df).$$

So as a result we have

$$\gamma \varepsilon(\beta-1)(\alpha-1-(\beta-1))+\beta(\beta-\alpha)(d,f)+(\alpha-\beta)(a,df)$$
$$=(\alpha-\beta)\left(\gamma \varepsilon (\beta-1)-\beta (d,f)+(a,df)\right).$$

Obtaining
$$(d_\alpha , f_\alpha)=\frac{\gamma \varepsilon (\beta-1) - \beta(d,f)+(a,df)}{\alpha - \beta}.$$
Then since
$(d,f)=\gamma \varepsilon + (d_\alpha, f_\alpha) + (d_\beta, f_\beta),$
we get that 
$$(d_\beta, f_\beta)=-\frac{\gamma \varepsilon(\alpha-1)-\alpha(d,f)+(a,df)}{\alpha-\beta},$$
completing the calculation. $\square$

The pseudo composition case was classified using the so-called associatior relation. We will confirm that its appearance stems from the stricter fusion law and will show that this is the case no matter the eigenvalues.
\begin{lemma}\label{assoc}
If $A$ is a $\mathcal{J(\alpha,\beta)}$-axial algebra, that satisfies the $(\star)$ condition, $a \in A$ is a primitive axis and $d,f \in A$, such that $d=(d,a)a+d_\alpha+d_\beta$, $f=(f,a)a+f_\alpha+f_\beta$, where $(d,a)=\gamma,\ (f,a)=\varepsilon$ and $d_\alpha f_\alpha \in A_1(a)$, then
$$a(df)+d(af)+f(da)=(\alpha+2\beta)df+\gamma((1+\alpha)af-2\alpha\beta f)+\varepsilon((1+\alpha)ad-2\alpha\beta d)+$$
$$+\Bigl((1-3\alpha-2\beta-4\alpha\beta)\gamma\varepsilon + (\alpha-2\beta+1)\Psi(\beta)-\beta (1-\alpha)\Psi(\alpha)\Bigr)a,$$
where $\Psi(t)=\frac{\gamma \varepsilon (t-1) - t(d,f)+(a,df)}{\alpha - \beta}$.
\end{lemma}

\textit{Proof.} 
We will prove this by explicitly computing $a(df)+d(af)+f(da)$:
$$a(df)=a(\gamma \varepsilon a + \alpha (\gamma f_\alpha +\varepsilon d_\alpha)+ \beta (\gamma f_\beta + \varepsilon d_\beta) +d_\alpha f_\alpha + d_\alpha f_\beta + d_\beta f_\alpha+ d_\beta f_\beta).$$

Applying distributivity,
$$\gamma \varepsilon a+ \alpha^2 (\gamma f_\alpha +\varepsilon d_\alpha)+ \beta^2 (\gamma f_\beta + \varepsilon d_\beta)+\beta(d_\alpha f_\beta + d_\beta f_\alpha) +a(d_\alpha f_\alpha) + a(d_\beta f_\beta).$$

Another element is
$$d(af)=d(a(\varepsilon a+ f_\alpha+f_\beta))=(\gamma a+d_\alpha+d_\beta)(\varepsilon a+ \alpha f_\alpha + \beta f_\beta)=$$
$$=\gamma \varepsilon a + \gamma (\alpha^2 f_\alpha + \beta^2 f_\beta)+ \varepsilon (\alpha d_\alpha + \beta d_\beta) + (d_\alpha+d_\beta)(\alpha f_\alpha+ \beta f_\beta). $$

By symmetry we derive
$$f(da)=\gamma \varepsilon a + \varepsilon(\alpha^2 d_\alpha + \beta^2 d_\beta)+ \gamma (\alpha f_\alpha+ \beta f_\beta)+ (f_\alpha + f_\beta)(\alpha d_\alpha +\beta d_\beta).$$

Therefore $$a(df)+d(af)+f(da)=3\gamma \varepsilon a+ 2\alpha^2(\gamma f_\alpha +\varepsilon d_\alpha)+ 2\beta^2(\gamma f_\beta+\varepsilon d_\beta) +\alpha (\gamma f_\alpha+ \varepsilon d_\alpha)+ \beta (\gamma f_\beta + \varepsilon d_\beta)+$$
$$+2\alpha f_\alpha d_\alpha +2 \beta f_\beta d_\beta + (\alpha+2\beta) (d_\alpha f_\beta+ d_\beta f_\alpha)+ a(d_\alpha f_\alpha)+ a(d_\beta f_\beta). $$

Then it becomes somewhat clear what needs to be subtracted. So we compute
$$a(df)+d(af)+f(da)-(\alpha+2\beta)df=(3-\alpha-2\beta)\gamma \varepsilon a +(\alpha^2-2\alpha \beta)(\gamma f_\alpha +\varepsilon d_\alpha)-\alpha \beta (\gamma f_\beta + \varepsilon d_\beta) + $$ 
$$+\alpha (\gamma f_\alpha+ \varepsilon d_\alpha)+\beta (\gamma f_\beta + \varepsilon d_\beta)+ a(d_\alpha f_\alpha)+(\alpha-2\beta)f_\alpha d_\alpha + a(f_\beta d_\beta)-\alpha (f_\beta d_\beta).$$ 

Note that $$(1+\alpha)ad-2\alpha\beta d= (1+\alpha)\gamma a + \alpha(1+\alpha)d_\alpha + \beta(1+\alpha)d_\beta - 2\alpha\beta \gamma a- \alpha 2\beta d_\alpha -\beta 2\alpha d_\beta =$$
$$=(1+\alpha-2\alpha\beta)\gamma a+\alpha (1+\alpha-2\beta)d_\alpha + \beta (1-\alpha)d_\beta.$$

Similarly, $(1+\alpha)af-2\alpha\beta f=(1+\alpha-2\alpha\beta)\varepsilon a+\alpha (1+\alpha-2\beta)f_\alpha + \beta (1-\alpha)f_\beta.$

Therefore 
$$a(df)+d(af)+f(da)-(\alpha+2\beta)df-\gamma((1+\alpha)af-2\alpha\beta f)-\varepsilon((1+\alpha)ad-2\alpha\beta d)=$$
$$=(1-3\alpha-2\beta-4\alpha\beta)\gamma\varepsilon a + a(d_\alpha f_\alpha)+(\alpha-2\beta)f_\alpha d_\alpha + a(f_\beta d_\beta)-\alpha(f_\beta d_\beta).$$

Since $d_\beta f_\beta \in A_1(a) \oplus A_\alpha(a)$, we have that
$$a(d_\beta f_\beta)-\alpha(d_\beta f_\beta)=-\beta (1-\alpha)\Psi(\alpha)a.$$
This transforms our equality into 
$$a(df)+d(af)+f(da)-(\alpha+2\beta)df-\gamma((1+\alpha)af-2\alpha\beta f)-\varepsilon((1+\alpha)ad-2\alpha\beta d)=$$
$$=(1-3\alpha-2\beta-4\alpha\beta)\gamma\varepsilon a + a(d_\alpha f_\alpha)+(\alpha-2\beta)f_\alpha d_\alpha -\beta (1-\alpha)\Psi(\alpha)a.$$
Finally, note that if $d_\alpha f_\alpha \in \langle a \rangle$, then 
$a(d_\alpha f_\alpha)=d_\alpha f_\alpha$ and $d_\alpha f_\alpha = (d_\alpha f_\alpha, a) a=\Psi(\beta)a$. $\square$

This demonstrates an effect, initially observed for the case of $\mathcal{PC}(\eta)$ algebras \cite{GorMamStar}. Its importance can be seen for the case when algebras are generated by $3$ or more primitive axes, since it imposes a relation on the products of generating elements, thus limiting the algebra's dimension (See Corollary \ref{pseudocompCor} below).

Another statement of interest is the Seress identity, which was used to great success in the classification of $3$-generated algebras of Jordan type by Gorshkov and Staroletov \cite{GorStarJordan}. However it is also not easily generalisable to our setting:

\begin{statement}[Seress Identity]\label{Seress}
Let $a$ be a primitive axis of a $\mathcal{J}(\alpha,\beta)$ algebra $A$, while $y \in A_1(a) + A_\alpha(a)$. Then for any $x \in A$ the identity
$$a(xy)=(ax)y$$ holds iff $\alpha=0$.
\end{statement} 

\textit{Proof.}
Assuming the identity holds, we see that for $0 \neq y \in A_\alpha(a)$ and $0\neq x \in A_1(a)$.
Then $$a(xy)-(ax)y=a(\lambda a y)-\lambda ay=\lambda(\alpha^2-\alpha)y=\lambda \alpha(\alpha-1)y.$$
implying that $\alpha=0$. 

If $\alpha=0$ this can be proved similarly to the proof given in \cite[Lemma 4.3]{Jortype} with a small correction, since the fusion law in this work is stricter. 

It is clear that if $y \in A_1(a)$ then the identity holds since
$$a(xy)=a(\lambda ax)=\lambda a (ax)=y(ax),$$
so by linearity we can assume $y \in A_0(a)$. Also by linearity, we can assume that $x \in A_t(a)$ for $t \in \{1, 0, \beta\}$. If $t=1$, then 
$$a(xy)=a(\lambda a y)=0=\lambda ay (a\lambda a)y=(ax)y.$$
If $t = 0$ then $xy \in A_1(a) +A_0(a)$ and we have that
$$a(xy)=a((a,xy)a+(xy)_0)=(a,xy)a=(ay,x)a=0=(ax)y,$$
while for $t=\beta$ $a(xy)=\beta(xy)=(ax)y$, since $xy \in A_\beta(a)$.
$\square$ 

From this we infer that the classification of finitely-generated $\mathcal{J}(\alpha,\beta)$-algebras requires some more general methods.

Statement \ref{Seress} also shows that the class of algebras that satisfy the Seress identity might be slightly greater than the class of algebras of Jordan type, since we have shown that said identity holds for a more general fusion law. Hence it may be interesting to find an algebra with such an identity, that is neither a Jordan algebra nor a Matsuo algebra.

\section*{Algebras, generated by $3$ axes}
The goal of this section is to study $3$-generated axial algebras of type $\mathcal{J}(\alpha, \beta)$. We will be using a Miyamoto-involution based approach, which is inspired by \cite{DeMetsEtal}.
Following this paper, we denote $[a]x=\tau_a(x)$, while $[a_1,\ldots a_l]b=[a_1]([a_2](\ldots [a_l]b))$. Note that this is equal to the action of composition of operators $a_1 a_2 \ldots a_l$, which is known to be associative.

We set out to prove the following theorem.

\begin{theorem*}
Assume that $\mathbb{F}$ is a field with $\alpha, \beta \in \mathbb{F}$ such that $\alpha \neq \beta$ and both of them are not equal to $1$. Then let $A$ be a primitive $\mathcal{J}(\alpha,\beta)$-algebra over $\mathbb{F}$, satisfying the $(\star)$ condition. If $A$ is generated by $3$ primitive axes $a,b,c$,  then $A$ is spanned by $a,\ b,\ c,\ [a]b,\ [a]c,\ [b]c,\ [a,b]c, [b,a]c, [c,a]b$. 

Moreover, if we denote $\delta=(a,b),\ \theta=(a,c),\ 
\gamma=(b,c)\ \omega=(a,[b]c)$, then the following is true for $a$:
$$A_\alpha(a)=\langle -2\delta a+b+[a]b, -2\theta a+c+[a]c,\ -2\omega a+[b]c+[a,b]c,\ \frac{t_1a-t(b+c)}{(\alpha-\beta)^2}+[b,a]c+[c,a]b\rangle,$$
$$A_\beta(a)=\langle -b+[a]b,\ -c+[a]c,\ -[b]c+[a,b]c,\ \frac{t}{(\alpha-\beta)^2}(b-c)-[b,a]c+[c,a]b\rangle,$$
where 
$t_1=-2 \beta^2 (-\omega + \delta + 2 \theta) - 4\beta \delta (-2 \theta + \gamma) + \alpha^2 (\omega - \delta (1 + 2 \gamma)) + 2 \alpha \delta \gamma + 2\alpha \beta (\theta + 2 \delta \theta - (2 + \delta) \gamma)$ and
$t=\alpha^2 (2 \omega + 4 \delta (\theta - \gamma) + 2 \theta - 2 \gamma + 1) - 2 \alpha \bigl(\beta (\omega - \theta + \gamma) + \omega + (4 \delta + 1) (\theta - \gamma)\bigr) - \beta^2 + 2 \beta (\omega - \theta + \gamma) + 4 \delta (\theta - \gamma)$.

\end{theorem*}

For the rest of this section by $A$ we denote an axial algebra of type $\mathcal{J}(\alpha,\beta)$, that satisfies the $(\star)$ condition and is generated by primitive axes $a,b,c$.

Our first order of business is to  obtain an upper bound on the dimension of $A$.

Consider the axial algebra $A$, generated by a finite set of axes denoted $S$.
We define $S\left[i \right]$ as $\langle \tau_{a_1}\tau_{a_2}\tau_{a_3}\ldots \tau_{a_l}(b)| l \leq i,\ a_i,b \in S\rangle $.

We will start by proving the following result for our fusion law:

\begin{prop}[c.f. \cite{DeMetsEtal}] \label{Sn}
If $S[n]=S[n+1]$ for some $n \geq 1$ then $A=S[n]$.
\end{prop}
\textit{Proof.}
Define the closure $C$ of $S=\{a,b,c\}$ to be the smallest subset of axes of $A$, containing $S$ such that for any $x \in C$ $\tau_x(C) \subseteq C$.

It can be shown that $A$ is spanned by $C$. Indeed, note that for any two axes $c,d \in A$ the subalgebra generated by these axes is spanned by $c,\ d,\ \tau_c(d)$. This means that $\langle C \rangle$ is multiplicatively closed in $A$, since the product of any two elements of $C$ lies in the subalgebra generated by these two elements and therefore lies in $ \langle C \rangle$. Since $C$ contains the generators of $A$ we conclude that $\langle C \rangle =A$.
By \cite[Lemma 3.5]{KhasMcShpec}, we have that $C=\{\tau_{a_1}\tau_{a_2}\ldots \tau_{a_l}(b),\ a_i,\ b \in S\}$. Also note that $S[n]=S[n+1]$ implies $S[n]=S[l]$ for all $l \geq n$ and this means that $\langle C \rangle =S[n]$, proving the proposition. $\square$

For further work we will compute $\tau_a(x)$ for a primitive axis $a$ and arbitrary $x$. This is equal to
$$\tau_a(x)=(a,x)a+x_\alpha - x_\beta.$$

We can derive from previous calculations that
$$x_\alpha=\frac{(a,x)(\beta-1)a-\beta x+ ax}{(\alpha-\beta)}, x_\beta=\frac{-(a,x)(\alpha-1)a+\alpha x-ax}{(\alpha-\beta)}.$$

Subtracting these and adding the third element yields us

$$\tau_a(x)=\frac{2(a,x)(\alpha-1)a-(\alpha+\beta)x+2ax}{\alpha-\beta}.$$

In particular, we have an interesting equality for axes $a$ and $b$

$$\tau_a(b)-\tau_b(a)=-\frac{2(a,b)(\alpha-1)+(\alpha+\beta)}{\alpha-\beta}(b-a).$$

Following \cite{DeMetsEtal} we will denote the obtained coefficient by $\varepsilon_{a,b}$.

$$\varepsilon_{a,b}=-\frac{2(a,b)(\alpha-1)+(\alpha+\beta)}{\alpha-\beta}.$$

As an example of usage of the square bracket notation the last identity can be rewritten as, $[a]b-[b]a=\varepsilon_{a,b}(b-a)$.

Another result known for axial algebras is the following: 

\begin{lemma}[cf Lemma 5.1 of \cite{Jortype}]
 
Let $\varphi \in Aut(A)$, then for an axis $a$ and $x \in A$
$$\tau_a^\varphi(x)=\tau_{\varphi(a)}(x),$$
where the left hand side is conjugated by $\varphi$.
\end{lemma}

\textit{Proof.}
Write $x$ as $\alpha a+x_\alpha+x_\beta$. Then since $\varphi$ is an automorphism, we have $\varphi(x)= \alpha \varphi(a)+ \varphi(x_\alpha)+\varphi(x_\beta)$, which is an expression of $x$ with respect to $\varphi(a)$, since $A_t(\varphi(a))=\varphi(A_t(a))$. 

Thus we need to show that the left hand side is equal to $\alpha \varphi(a)+\varphi(x_\alpha)-\varphi(x_\beta)$. Now, the left hand side is $\varphi \circ \tau_a \circ \varphi^{-1}$, which can be rewritten into $\alpha \varphi(\tau_a(a))+\varphi(\tau_a(x_\alpha))+\varphi(\tau_a(x_\beta))$, which is precisely the right hand side. $\square.$

From this general result we get that
$$\left[ a,b, a \right]=\left[\tau_a(b)\right].$$
This means that we only need to work with Miyamoto involutions that are induced by generating axes $a,\ b,\ c$.

We will be constructing a spanning set for $A$ by constructing spanning sets for $S[i]$ until we reach $i$ such that $S[i]=S[i+1]$ which, per Proposition \ref{Sn}, will show that we are done.
We start with $S[0]$, with its spanning set consisting of $3$ elements.
Then we have $S[1]$ which due to the equality $[a]b-[b]a=\varepsilon_{a,b}(b-a)$ presented above adds, for instance $[a]b, [a]c, [b]c$ to the global spanning set. 

$S[2]$ in turn adds only $[a,b]c,\ [c,a]b,\ [b,a]c$ due to the following identity

$$[a,b]a=\varepsilon_{a,b}a+b-\varepsilon_{a,b}[a]b.$$ 

This holds since
$$[a,b]a=[a]([b]a-[a]b+[a]b)=[a](\varepsilon_{a,b}(a-b)+[a]b).$$
Additionally,
$$[a,b]c=[a][b]c=[a]([b]c-[c]b+[c]b)=[a](\varepsilon_{b,c}(c-b)+[c]b)=\varepsilon_{c,b}[a](c-b)+[a,c]b, \eqno{(**)}$$
implying $[a,b]c$ and $[a,c]b$ are linearly dependent in $S[2]$ (recall that $S[2]$
contains $S[1]$ and $S[0]$).
Now we will show that $S[2]$ is invariant under the actions of Miyamoto involutions.

Consider the word $[x,y,z]w$ where $x,\ y,\ z,\ w \in \{a,\ b,\ c\}$. To avoid the cases, where this word is trivially in $S[2]$, no two consecutive letters must be equal. This leaves us with the following words (up to the permutation of letters)

$$[a,b,a]c,\ [a,b,a]b,\ [a,b,c]b,\ [a,b,c]a.$$
Words $[a,b,a]b,\ [a,b,c]b$ can be shown to be linearly dependent with $[a,b,b]a$ and $[a,b,b]c$ respectively via a trick presented above at $(**)$.

Word $[a,b,a]c$ can be shown to be equal to $\varepsilon_{\tau_{a}(b),c}(c-[a]b)+[c,a]b$ and the word of form $[a,b,c]a$ is equal to $[a,b,a]c+\varepsilon_{a,c}[a,b](a-c)$ \cite[Lemma 3.11]{DeMetsEtal}.

This shows us that $S[3]=S[2]$, confirming the set of axes from the Theorem indeed spans $A$. $\square$

During the course of the proof we have obtained a spanning set for the algebra $A$, which consists of primitive axes.
$$a,\ b,\ c,\ [a]b,[a]c, [b]c, [a,b]c,\ [b,a]c,\ [c,a]b.$$

From this and the expression for Miyamoto involutions we get the following 
\begin{coro}
 Let $A$ be a primitive $\mathcal{J}(\alpha,\beta)$-axial algebra satisfying the $(\star)$ condition. If $A$ is generated by $3$ primitive axes $a,\ b,\ c$, then $A$ is spanned by
$a,\ b,\ c,\ ab,\ ac,\ bc,$ $a(bc),\ b(ac),\ c(ab)$.
\end{coro}

\begin{coro}\label{pseudocompCor}
    If we assume that for a primitive axis $a$ $A^2_\alpha(a) \subseteq A_1(a)$, then the algebra is spanned by $a,\ b,\ c,\ ab,\ ac,\ bc,\ a(bc),\ b(ac)$.
\end{coro}
\textit{Proof.} Note that by Lemma \ref{assoc} we have that $a(bc)+b(ac)+c(ab)$ can be expressed as a linear combination the other elements of the above basis. $\square$

This corollary holds for any generating primitive axis (assuming the corresponding basis elements are replaced).

We can now use the newly obtained dimension bound to study the algebra in more detail. We will start with the computation of the Gram matrix.

It is greatly simplified with the following

\begin{lemma}
Let $a$ be a primitive axis in $A$.
Then $([a]x,y)=(x,[a]y)$ for any $x,y \in A$.
\end{lemma}
\textit{Proof.} It is clear that
$$([a]x,y)=((a,x)a+x_\alpha-x_\beta,(a,y)a+y_\alpha+y_\beta)=(a,x)(a,y)+(x_\alpha,y_\alpha)-(x_\beta,y_\beta). \; \square$$

This result makes it incredibly easy to work with the basis we have constructed.

\begin{prop}
Let $A$ be a $\mathcal{J}(\alpha,\beta)$-algebra, satisfying the $(\star)$ condition, that is generated by primitive axes $a,b,c$. Denote the values of the Frobenius form as follows:
$\delta=(a,b),\ \theta=(a,c),\ \gamma=(b,c),\ \omega=(a,[b]c)$. Then the Gram matrix of $A$ is presented in the Table 2.
\end{prop}

\begin{sidewaystable}[h]
    \centering
\scalebox{0.78}{
\begin{tabular}{|C||C|C|C|C|C|C|C|C|C|}
\hline
(,) & a & b & c& [a]b& [a]c &[b]c&[a,b]c & [b,a]c & [c,a]b \\\hline
a & 1 & \delta & \theta & \delta & \theta & \omega & \omega & \varepsilon_{a,b}\theta+\gamma- &  \varepsilon_{a,c}\delta+\gamma-\\
& & & & & & & & -\varepsilon_{a,b}(\omega-\varepsilon_{a,b}(\theta-\gamma))& -\varepsilon_{a,c}(\omega-\varepsilon_{a,b}(\theta-\gamma))\\ \hline
b & & 1 & \gamma & \delta-\varepsilon_{a,b}(\delta-1) & \omega-\varepsilon_{a,b}(\theta-\gamma) & \gamma &  \theta-\varepsilon_{a,b}(\omega-\gamma) & \omega-\varepsilon_{a,b}(\theta-\gamma) & \varepsilon_{a,b}\gamma+\theta-\varepsilon_{a,b}\omega-\\
& & & & & & & & &-\varepsilon_{b,c}((\omega-\varepsilon_{a,b}(\theta-\gamma)) \\
& & & & & & & & & -(\delta-\varepsilon_{a,b}(\delta-1)))\\ \hline
c & & & 1 & \omega-\varepsilon_{a,b}(\theta-\gamma) & \theta-\varepsilon_{a,c}(\theta-1) & \gamma-\varepsilon_{b,c}(\gamma-1)& \varepsilon_{b,c}\theta+\delta- & \varepsilon_{b,c}\theta+\delta- & \omega-\varepsilon_{a,b}(\theta-\gamma) \\
& & & & & & & -\varepsilon_{b,c}(\omega-\varepsilon_{b,c}(\theta-\delta))-&-\varepsilon_{b,c}(\omega-\varepsilon_{b,c}(\theta-\delta))- & \\
& & & & & & & -\varepsilon_{a,c}(\omega-(\gamma-\varepsilon_{b,c}(\gamma-1)))& -\varepsilon_{a,c}(\omega-(\gamma-\varepsilon_{b,c}(\gamma-1))) & \\ \hline
[a]b & & & & 1 & \gamma & \theta-\varepsilon_{a,b}(\omega-\gamma) & \gamma & \theta+\varepsilon_{a,b}(\omega-\varepsilon_{a,b}(\theta-\gamma))-& \varepsilon_{\tau_a(c),b}(1-(\omega-\varepsilon_{a,b}(\theta-\gamma)))- \\
& & & & & & & & -\varepsilon_{a,b}(\varepsilon_{a,b}\theta+\gamma) &-\varepsilon_{a,b}(\theta-\gamma)+\omega \\
& & & & & & & & -\varepsilon^2_{a,b}(\omega-\varepsilon_{a,b}(\theta-\gamma))& \\ \hline
[a]c & & & & & 1 &\varepsilon_{b,c}\theta+\delta- &\gamma-\varepsilon_{b,c}(\gamma-1) & (1-\varepsilon_{\tau_a(b),c})(\omega-\varepsilon_{a,b}(\theta-\gamma))+ &\varepsilon_{\tau_{a}(c),b}(\gamma-(\theta-\varepsilon_{a,c}(\theta-1)))+\\ 
& & & & & &-\varepsilon_{b,c}(\omega-\varepsilon_{b,c}(\theta-\delta))-& &+\varepsilon_{\tau_a(b),c}  &+([a]c,[b]c) \\
& & & & & &  -\varepsilon_{a,c}(\omega-(\gamma-\varepsilon_{b,c}(\gamma-1)))& & & \\ \hline
[b]c & & & & & & 1& \varepsilon_{\tau_b(a),c} (1-\omega)+\omega & \theta-\varepsilon_{a,c}(\theta-1) &\varepsilon_{b,c}(\omega-\varepsilon_{a,b}(\theta-\gamma))+ \delta-\\
& & & & & & & & &-\varepsilon_{a,b}(\delta-1)-\varepsilon_{b,c}(b,[c,a]b))\\ \hline
[a,b]c& & & & & & & 1&\varepsilon_{\tau_a(b),c}(\gamma-\theta-\varepsilon_{b,c}(\gamma-1)+ &\varepsilon_{\tau_a(c),b}(\gamma-([b]c,[a]c))+ \\ 
& & & & & & & &+\varepsilon_{a,b}(\omega-\gamma))+([b]c,[c,a]b) & +\theta-\varepsilon_{a,c}(\theta-1)\\\hline
[b,a]c & & & & & & & &1 &\varepsilon_{\tau_b(c),a}(\theta-([a]c,[b]c))+\gamma- \\ 
& & & & & & & & &-\varepsilon_{b,c}(\gamma-1)+\\ 
& & & & & & & & &+\varepsilon_{a,b}\varepsilon_{b,c}(\omega-\varepsilon_{a,b}(\theta-\gamma)-\\
& & & & & & & & &-([a]c,b[c]))+\theta-\varepsilon_{a,c}(\theta-1)-\\
& & & & & & & & &-([b,a]c,[a]c)-\\
& & & & & & & & &-\varepsilon_{a,c}(([a]c,c)-\theta)\\ \hline
[c,a]b & & & & & & & & &1 \\ \hline

\end{tabular}}
 \caption{The Gram Matrix values for a primitive $\mathcal{J}(\alpha,\beta)$-algebra generated by axes $a,\ b,\ c$.}
\end{sidewaystable}
The table is self-referential in order to preserve readability.

Since our basis consists of primitive axes we can also construct its multiplication table, similar to the $2$-generated case presented above. This can be done as follows: take any two of the basis elements, which we will denote by $d$ and $f$. Being axes, they form an at most $3$-dimensional axial algebra, which is known to be a span of $d,f,$ and $\tau_d(f)$. We can then use the rewriting rules proposed in  \cite{DeMetsEtal} in order to transform $\tau_d(f)$ and find $df$ via the expressing $\tau_d(f)$ as a combination of $d,f, df$.

For instance, we can find $[a]b[a]c$ via $[[a]b,a]c=[[a]b][a]c$, which is known to be equal to 
$$[a,b,a,a]c=[a,b]c,$$ thus getting that 
$$[a]b[a]c=\frac{1}{2}(\alpha-\beta)[a,b]c-2\gamma(\alpha-1)a[b]+(\alpha+\beta)[a]c,$$
where $\gamma=(b,c)$.

This example is rather simple so we will show another example, involving a more complicated derivation.

The product $[a]b [b,a]c$ can be expressed via $[a]b,\ [b,a]c$ and
$$[a,b,a,b,a]c=[a,b](\varepsilon_{\tau_a(b),c}(c-[a]b)+[c,a]b)=$$
$$=[a]\varepsilon_{\tau_a(b),c}\Bigl([b]c-\varepsilon_{a,b}(b-[b]a)+a\Bigr)+$$
$$+\varepsilon_{\tau_b(c),a}(a-[b]c)+[a,b]c+\varepsilon_{a,b}(\varepsilon_{b,c}(b-[b]c)+c-[b,c]a)$$
$$=\varepsilon_{\tau_a(b),c}\Bigl([a,b]c-\varepsilon_{a,b}([a]b-(\varepsilon_{a,b}(a-[a]b)+b))+a\Bigr)+\varepsilon_{\tau_b(c),a}(a-[a,b]c)+[b]c+$$
$$+\varepsilon_{a,b}\Bigl(\varepsilon_{b,c}([a]b-[a,b]c)+[a]c-\bigl(\varepsilon_{\tau_a(b),c}(c-[a]b)+[c,a]b+\varepsilon_{a,c}(\varepsilon_{a,b}(a-[a]b)+b-[a,b]c)\bigr)\Bigr).$$

Which, after collecting the terms leaves us with

$$[a,b,a,b,a]c=(\varepsilon_{\tau_a(b),c}(1-\varepsilon_{a,b})+\varepsilon_{\tau_b(c),a}-\varepsilon^2_{a,b}\varepsilon_{a,c})a+(\varepsilon_{a,b}(\varepsilon_{a,c}-\varepsilon_{\tau_a(b),c}))b-(\varepsilon_{a,b}\varepsilon_{\tau_{a}(b),c})c+$$
$$+\Bigl(\varepsilon_{a,b}(\varepsilon_{\tau_a(b),c}(\varepsilon_{a,b}-1)+\varepsilon_{b,c}+\varepsilon_{a,c})-\varepsilon_{\tau_a(b),c}\Bigr)[a]b+\varepsilon_{a,b}[a]c+[b]c+$$
$$+(\varepsilon_{\tau_a(b),c}-\varepsilon_{\tau_b(c),a}-\varepsilon_{a,c}(\varepsilon_{a,b}+1))[a,b]c-\varepsilon_{a,b}[c,a]b.$$

Then, substituting this into the formula yields us 
$$[a][b,a]c=\frac{1}{2}(\alpha-\beta)[a,b,a,b,a]c-([a]b,[b,a]c)(\alpha-1)a[b]-\frac{(\alpha+\beta)}{2}[b,a]c.$$

As it can be seen these derivations involve a great number of steps and require a certain amount of caution. Further on, said derivations become even longer and do not require usage of any new techniques and, hoping that Proposition 7 and the examples provided above are enough to prove and demonstrate the effectiveness of the algorithm, we will refrain from presenting all of the computations in this paper. However both the multiplication table and the Gram matrix can be found at \cite{git}, reproduced in the computer algebra system GAP \cite{GAP}.

We will also present some of the simpler products in order to prove

\begin{prop}
    Let $A$ be a $\mathcal{J}(\alpha,\beta)$-algebra, satisfying the $(\star)$ condition, that is generated by primitive axes $a,b,c$. Then if we denote $\delta=(a,b),\ \theta=(a,c),\ \gamma=(b,c)\ \omega=(a,[b]c)$, then the following is true for a primitive axis $a$:
$$A_\alpha(a)=\langle -2\delta a+b+[a]b, -2\theta a+c+[a]c,\ -2\omega a+[b]c+[a,b]c,\ \frac{t_1a-t(b+c)}{(\alpha-\beta)^2}+[b,a]c+[c,a]b\rangle,$$
$$A_\beta(a)=\langle -b+[a]b,\ -c+[a]c,\ -[b]c+[a,b]c,\ \frac{t}{(\alpha-\beta)^2}(b-c)-[b,a]c+[c,a]b\rangle,$$
where 
$t_1=-2 \beta^2 (-\omega + \delta + 2 \theta) - 4\beta \delta (-2 \theta + \gamma) + \alpha^2 (\omega - \delta (1 + 2 \gamma)) + 2 \alpha \delta \gamma + 2\alpha \beta (\theta + 2 \delta  \theta  - (2 + \delta) \gamma)$ and
$t=\alpha^2 (2 \omega + 4 \delta (\theta - \gamma) + 2 \theta - 2 \gamma + 1) - 2 \alpha \bigl(\beta (\omega - \theta + \gamma) + \omega + (4 \delta + 1) (\theta - \gamma)\bigr) - \beta^2 + 2 \beta (\omega - \theta + \gamma) + 4 \delta (\theta - \gamma)$.
\end{prop}

For this it suffices to find the eigenvectors of the matrix, consisting of product vectors, which will look as follows:
$$ab=\frac{1}{2}(\alpha-\beta)[a]b-(\alpha-1)\delta a+\frac{(\alpha+\beta)}{2}b,\quad ac=\frac{1}{2}(\alpha-\beta)[a]c-(\alpha-1)\theta a+\frac{(\alpha+\beta)}{2}c,$$
$$a[a]b=\frac{1}{2}(\alpha-\beta)b-(\alpha-1)\delta a+\frac{(\alpha+\beta)}{2}[a]b, \quad a[a]c=\frac{1}{2}(\alpha-\beta)c-(\alpha-1)\theta a+\frac{(\alpha+\beta)}{2}[a]c,$$
$$a[b]c=\frac{1}{2}(\alpha-\beta)[a,b]c-(\alpha-1)\omega a+\frac{(\alpha+\beta)}{2}[b]c,\; \; a[a,b]c=\frac{1}{2}(\alpha-\beta)[b]c-(\alpha-1)\omega a+\frac{(\alpha+\beta)}{2}[a,b]c,$$
$$a[b,a]c=\frac{1}{2}(\alpha-\beta)(\varepsilon_{\tau_{a}(b),c}(c-[a]b)+[c,a]b)-(\alpha-1)(\varepsilon_{a,b}\theta +\gamma -\varepsilon_{a,b}(\omega-\varepsilon_{a,b}(\theta-\gamma)))a+\frac{(\alpha+\beta)}{2}[b,a]c,$$
$$a[c,a]b=\frac{1}{2}(\alpha-\beta)(\varepsilon_{\tau_{a}(c),b}(b-[a]c)+[b,a]c)-(\alpha-1)(\varepsilon_{a,c}\delta+\gamma-\varepsilon_{a,c}(\omega-\varepsilon_{a,b}(\theta-\gamma)))a+\frac{(\alpha+\beta)}{2}[c,a]b.$$

The eigenvectors can then be found via any suitable computer algebra system.

We note that this method for constructing the multiplication table is applicable for any finitely-generated $\mathcal{J}(\alpha,\beta)$ algebras as long as they possess a basis consisting of axes and a sufficient number of rewriting identities. Due to this and the fact that \cite{DeMetsEtal} contains rewriting rules for the $4$-generated algebras of Jordan type we state the following

\begin{conj}
    Let $A$ be a $\mathcal{J}(\alpha,\beta)$-algebra, generated by $4$ primitive axes. Then $dim(A)\leq 81$.
\end{conj}

This, if proven true, may suggest a broader hypothesis that the algebras of Jordan type serve as an upper bound for $\mathcal{J}(\alpha,\beta)$-algebras dimension-wise.

\end{document}